\documentclass[12pt,reqno]{amsart}
\usepackage[a4paper,margin=1in]{geometry}
\usepackage[T1]{fontenc}
\usepackage{lmodern}
\usepackage{microtype}
\usepackage{amsmath,amssymb,amsthm,mathtools}
\usepackage{xcolor}
\usepackage{hyperref}
\usepackage[nameinlink,capitalise,noabbrev]{cleveref}
\usepackage{fancyhdr}

\hypersetup{
  colorlinks=true,
  linkcolor=blue,
  citecolor=green!35!black,
  urlcolor=blue!65!black,
}

\newtheorem{theorem}{Theorem}[section]
\newtheorem{proposition}[theorem]{Proposition}
\newtheorem{lemma}[theorem]{Lemma}
\theoremstyle{definition}
\newtheorem{definition}[theorem]{Definition}
\numberwithin{equation}{section}
\newtheorem{Remark}[theorem]{Remark}

\newcommand{\Ric}{\operatorname{Ric}}

\newcommand{\dd}{\,\mathrm d}
\newcommand{\IR}{\mathbb{R}}
\newcommand{\Sn}{\mathbb S}
\newcommand{\abs}[1]{\left\lvert #1\right\rvert}
\newcommand{\inner}[2]{\left\langle #1,#2\right\rangle}
\newcommand{\e}{\mathrm e}
\newcommand{\Oterm}{\mathrm O}

\title{
Uniqueness of the Bryant Soliton via the Positive Mass Theorem
}

\author[P.-Y.~Chan]{Pak-Yeung Chan}
\address[]{Department of Mathematics, National Tsing Hua University, Hsin-Chu, Taiwan}
\email{pychan@math.nthu.edu.tw}

 \author[Y.~Lai]{Yi Lai}
\address[]{Department of Mathematics, Rowland Hall, University of California Irvine, CA 92617, USA
}
\email{ylai25@uci.edu}

\begin{document}
\begin{abstract}
For all $n\ge 4$,
we prove that every $n$-dimensional complete noncompact
steady gradient Ricci soliton asymptotic to the cylinder $\mathbb R\times\mathbb S^{n-1}$ is, up to scaling, isometric to the Bryant soliton.

Our proof combines the exponential Bryant asymptotics of Bamler--Lai~\cite{BamlerLai}
with the positive mass theorem in all dimensions through a scalar curvature comparison. 
\end{abstract}

\maketitle

\section{Introduction}
A steady gradient Ricci soliton is a connected complete Riemannian
manifold $(M^n,g)$ together with a function $f\in C^\infty(M)$ satisfying
\begin{equation}\label{eq:steady-soliton}
 \Ric_g=\nabla_g^2f.
\end{equation}
For every $n\ge3$, the Bryant soliton~\cite{Bryant05} is a complete,
rotationally symmetric, non-collapsed steady gradient Ricci soliton on
$\mathbb R^n$ with positive curvature operator. It is asymptotic to the
cylinder $\mathbb R\times\mathbb S^{n-1}$ (see \cref{defn:asym_round}).
This raises the following question: Does asymptotic cylindricity uniquely
determine the Bryant soliton, up to scaling, among complete noncompact
steady gradient Ricci solitons?
(See, e.g., \cite[\S~5.2]{HaslhoferICM2026} and
\cite[Remark~4.13]{ZhaoZhuSurvey2026}.)

In dimension three, Brendle~\cite{Brendle2013} answered this question
affirmatively. His theorem settled the steady soliton case of
Perelman's conjecture, posed in the work leading to his proof of
the Poincar\'e conjecture by Ricci flow with
surgery~\cite{Perelman2002, Perelman2003,PerelmanFiniteExtinction}: up to scaling, the
Bryant soliton is the only noncompact three-dimensional non-collapsed
ancient solution with bounded positive curvature~\cite[Remark~11.9]{Perelman2002}.

Brendle later proved the full conjecture by establishing rotational
symmetry~\cite{Brendle2020}. Bamler--Kleiner gave an alternative proof of
this symmetry result~\cite{BamlerKleiner}.

In dimensions $n\ge4$, this question is more delicate because
asymptotic cylindricity does not directly imply nonnegative sectional
curvature, which is automatic for complete ancient Ricci flows in
dimension three~\cite[Corollary~2.4]{Chen2009}.
Under the additional assumption of positive sectional curvature,
Brendle~\cite{Brendle2014} proved uniqueness. In this paper, we answer the
question affirmatively in all dimensions $n\ge4$, without assuming any curvature sign or
prescribing the topology of $M$.

\begin{theorem}[Uniqueness of the Bryant soliton]\label{thm:main}
Let $n\ge4$. Every $n$-dimensional complete noncompact steady gradient
Ricci soliton asymptotic to the cylinder
$\mathbb R\times\mathbb S^{n-1}$ in the sense of \cref{defn:asym_round}
is, up to scaling, isometric to the Bryant soliton.
\end{theorem}

There are many more related rigidity results.
For example, Deng--Zhu~\cite{DengZhu2020,DengZhu2020a} proved uniqueness of
the Bryant soliton for non-collapsed steady solitons with linear scalar
curvature decay under suitable nonnegative curvature assumptions.
Zhao--Zhu~\cite{ZhaoZhu} replaced Brendle's curvature hypothesis by a weaker
global curvature-pinching condition, while Law~\cite{Law} obtained rigidity
under a bound on the dimension of the weighted $L^2$ kernel of the stability
operator.
See also~\cite{BamlerChowDengMaZhang, Brendle2011, BrendleNaff, ChanMaZhang, DengZhu2019,
ZhaoZhu2023, ZhaoZhu2024, ZhaoZhu2025,
CaoChen, CaoYu, CaoXie, CatinoMantegazza, Kim2017, Kim2025, Li2025,
DengZhu2021, Hebbar, RibeiroTran2026, DengWang2026}.

We note that the cylindrical hypothesis cannot be replaced by
non-collapsing, even for steady solitons with positive curvature
operator. In every dimension $n\ge4$, the flying wings of the second-named
author~\cite{Lai2024} give non-collapsed steady solitons on $\mathbb R^n$ distinct
from the Bryant soliton; see also \cite{ChanLaiLee}.

By \cref{thm:main}, asymptotic cylindricity determines the Bryant
soliton up to scaling. This raises the analogous uniqueness question for
quotient cylinders $\mathbb R\times(\mathbb S^{n-1}/\Gamma)$, where
$\Gamma$ is a finite group acting freely by isometries on $\mathbb S^{n-1}$.
For nontrivial $\Gamma$, the corresponding quotient of the Bryant soliton
is singular at the tip. In contrast, Appleton~\cite{Appleton} constructed
smooth, non-collapsed, non-K\"ahler steady gradient Ricci solitons on the
total spaces of the line bundles
$\mathcal O(-k)\to\mathbb{CP}^q$, where $q\ge1$ and $k\ge q+2$. These
solitons are asymptotic to the quotient cylinders
$\mathbb R\times(\mathbb S^{2q+1}/\mathbb Z_k)$. 
Recently, Ribeiro--Tran~\cite{RibeiroTran2026} proved uniqueness within
Appleton's ansatz in dimension four. Deng--Wang~\cite{DengWang2026} and
Law~\cite{Law2026} established uniqueness within the ansatz in every even
dimension \(n\ge4\).

 It remains
open whether every complete smooth steady gradient Ricci soliton with
such lens-space cylindrical asymptotics is, up to scaling, isometric
to the corresponding Appleton soliton.

\subsection*{Sketch of proof} 
Let $(M_B,g_B,f_B)$ denote the Bryant soliton metric and the potential, while $(M,g,f)$ denote that of the unknown soliton.
We use the scale normalization \eqref{eq:normalization} and fix the
additive constants by $\min_M f=\min_{M_B} f_B=0$.
Our starting point is Bamler--Lai's theorem on exponential Bryant
asymptotics~\cite[Theorem~1.3]{BamlerLai}. In \cref{s:exp} we verify that
the canonical Ricci flow of $g$ has dominant linear mode and satisfies the
remaining hypotheses. It follows that $g$ converges exponentially to $g_B$
at infinity. The soliton identities then give
exponential convergence of the potentials, up to an additive constant
(\cref{lem:potential-alignment}).
This gives the control at infinity needed for the scalar curvature comparison.

Our scalar curvature comparison is inspired by sliding arguments based on
the avoidance principle in mean curvature flow; see, for example,
Bamler--Lai's two-sided comparison with translates of the bowl
soliton~\cite[proof of Theorem~6.1]{BamlerLaiMCF1}.
In that proof, the asymptotic estimate ensures that any fixed nonzero
shift along the translating direction separates the two ends outside a compact
set. This rules out contact at infinity, so strong avoidance yields uniqueness.
However, the tensor maximum principle for Ricci flow does
not give an analogous comparison between the two soliton metrics.

Instead, we compare the scalar curvature of the unknown soliton with that
of the Bryant soliton. More precisely, for each shift $c\ge0$, we match
$p\in M$ with any point on the Bryant level set $\{f_B=f(p)+c\}$. The
scalar curvature difference is
\[
 R_g(p)-r_B(f(p)+c),
\]
where $r_B(T)$ denotes the Bryant scalar curvature on $\{f_B=T\}$. The exponential
asymptotics determine a constant $c_\infty\in\IR$ for which $R_g-r_B(f+c_\infty)$
is exponentially small. For fixed $c\ge0$ with $c>c_\infty$, the Bryant
asymptotics give
\[
 r_B(f+c_\infty)-r_B(f+c)\sim\frac{(n-1)(c-c_\infty)}{2f^2}
 \qquad\text{as }f\to\infty.
\]
This positive difference dominates the exponential error and gives
$R_g-r_B(f+c)>0$ near infinity. Starting from large $c$, for which this
inequality holds globally, we decrease $c$. By computation we can show
\[
 L_c\bigl(R_g-r_B(f+c)\bigr)\le0,
\]
where $L_c$ is a linear elliptic operator. At a first interior contact,
the difference is nonnegative and vanishes at an interior point. The
strong minimum principle for nonnegative supersolutions then gives
$R_g\equiv r_B(f+c)$. This is impossible for $c>c_\infty$ since the
difference is positive near infinity.
If $c_\infty<0$, the
equality $R_g=r_B(0)=1$ at a minimum of $f$ would give such a contact at
$c=0$. So
$c_\infty\ge0$ and $R_g\ge r_B(f+c_\infty)$.

The avoidance argument for the bowl soliton works in both directions,
whereas our sliding argument gives only a one-sided comparison.
To obtain equality, we use a conformal change and the positive mass theorem.
More precisely, we first
choose the function $u_B:[0,\infty)\to(0,\infty)$ so that
$(u_B^{\frac{4}{n-2}}\circ f_B)g_B$ is Euclidean.
Since $c_\infty\ge0$, the conformal factor $u_B(f+c_\infty)$
is defined on all of $M$.
For the complete metric $\widetilde g=u_B^{\frac{4}{n-2}}(f+c_\infty)g$, we have
\[
 u_B^{\frac{n+2}{n-2}}(f+c_\infty)R_{\widetilde g}
 =\bigl[R_g-r_B(f+c_\infty)\bigr]Q_n(f+c_\infty),
 \qquad Q_n>0.
\]
Since $R_g\ge r_B(f+c_\infty)$, this gives $R_{\widetilde g}\ge0$.
The same exponential asymptotics imply that $\widetilde g$ is asymptotically
flat and has zero ADM mass $m_{\mathrm{ADM}}(\widetilde g)=0$.
By the rigidity of the positive mass theorem by Schoen--Yau~\cite{SchoenYau79,SchoenYau81,Schoen89}
for $n\le7$ and Khuri--Wang--Wang~\cite{KhuriWangWang} in higher dimensions,
building on Brendle--Wang~\cite{BrendleWang}, $(M,\widetilde g)$ is isometric
to Euclidean space (see also \cite{SchoenYau79, SchoenYau81, Witten81, Bartnik, Lohkamp, SchoenYau22, CMSW25,
BHHSZ26} and references therein). Thus $g$ is locally conformally flat, and the
classification of Cao--Chen~\cite[Theorem~1.2]{CaoChen} and Catino--Mantegazza~\cite[Proposition~2.3]{CatinoMantegazza} identifies
$(M,g)$ with the Bryant soliton up to scaling.

\subsection*{Acknowledgment} 
P.-Y. Chan is  supported by the Yushan Young Fellow Program of the Ministry of Education (MOE), Taiwan (MOE-108-YSFMS-0004-012-P1), and by NSTC grants 113-2115-M-007 -014 -MY2 and 115-2628-M-007 -006. Y. Lai is supported by NSF grants DMS-2506832 and DMS-2543311 (CAREER).

\subsection*{Disclosure on AI assistance}

The authors used ChatGPT Pro to explore possible proof strategies and polish the exposition. It initially proposed a much longer and more involved argument. Inspired by the idea of using the positive mass theorem in that argument,\footnote{Note that the conformal change and positive mass theorem approach was also used by
Bunting and Masood-ul-Alam \cite{BuntingMasood} to give an alternative proof of uniqueness of
Schwarzschild spacetime \cite{Robinson1977}. Beig and Simon~\cite[\S~1]{BeigSimon} used a related
conformal comparison with a spherical reference solution to prove uniqueness for static perfect fluids.}
the authors found the much simpler proof presented here. The authors take full responsibility for the mathematical content. No AI-generated text appears in the final manuscript.

\section{Preliminaries}\label{s:prelim}
We begin by recalling some basic properties of gradient steady Ricci solitons and the definition of asymptotic cylindricality (\cref{defn:asym_round}). We then provide several equivalent characterizations of asymptotic cylindricality (\cref{lem:asym-round-equivalence}). Finally, we conclude this section by recalling the asymptotic behavior and ODEs satisfied by various geometric quantities of the Bryant steady soliton.

Hamilton~\cite{Hamilton95} showed that
\(R_g+\abs{\nabla_g f}_g^2\) is constant on every steady gradient Ricci
soliton. 
By \cite{Zhang2009}, the flow of \(-\nabla_g f\) is complete (see also \cite{Chen2009}).
Let \(\varphi_t\) denote this flow, with \(\varphi_0=\mathrm{id}\), and set
\begin{equation}\label{eq:canonical-flow}
 g(t):=\varphi_t^*g,
 \qquad
 \frac{\dd}{\dd t}\varphi_t=-\nabla f\circ\varphi_t.
\end{equation}
Then $g(t)$ satisfies the Ricci flow equation, i.e.
\(\partial_tg(t)=-2\Ric_{g(t)}\), and we call \(g(t)\) the
\textit{canonical Ricci flow} associated with \(g\).

The following steady-soliton identities are standard; see, for example,
\cite[Chapter~4]{ChowLuNi}:
\begin{align}
 \Delta f&=R_g,\label{eq:trace-soliton}\\
 \nabla R_g&=-2\Ric(\nabla f,\cdot),\label{eq:gradR}\\
 \Delta_fR_g&=-2\abs{\Ric}^2,\label{eq:driftR}
\end{align}
where \(\Delta_f:=\Delta+\nabla_{\nabla f}\).

By \cite[Corollary~2.5]{Chen2009}, \(R_g\ge0\). The strong minimum
principle applied to \eqref{eq:driftR} then gives \(R_g>0\) on a
non-Ricci-flat soliton.
The constant in Hamilton's
identity is therefore positive, and after rescaling we assume
\begin{equation}\label{eq:normalization}
 R_g+\abs{\nabla_g f}_g^2=1.
\end{equation}
In particular, \(\abs{\nabla_g f}_g\le1\). Moreover, combining \eqref{eq:trace-soliton} and \eqref{eq:normalization} we have $\Delta_f f=1$.

\medskip
We now define asymptotic cylindricity for a steady soliton.
See also~\cite{Brendle2013,Brendle2014,DengZhu2020,BamlerLai, BamlerChowDengMaZhang} for related
formulations of cylindrical asymptotics.
\begin{definition}[Asymptotic cylindricity]\label{defn:asym_round}
A complete connected noncompact non-Ricci-flat steady soliton \((M,g,f)\) of dimension $n\ge3$ is
\emph{asymptotically cylindrical}, or \emph{asymptotic to
\(\mathbb R\times\mathbb S^{n-1}\)}, if for every sequence $p_i\to\infty$,
with $r_i=(n-1)/(2R_g(p_i))$, we have
\[
 (M,r_i^{-1}g,p_i)\longrightarrow
 (\mathbb R\times\mathbb S^{n-1},\dd z^2+2(n-2)g_0,o)
\]
in the pointed smooth Cheeger--Gromov topology, where $o$ is a point on the cylinder and $g_0$ is
the standard metric of sectional curvature one on $\mathbb S^{n-1}$.
\end{definition}

The following lemma shows that \cref{defn:asym_round} is equivalent to two
standard formulations of asymptotic cylindricity for the canonical Ricci flow.

\begin{lemma}\label{lem:asym-round-equivalence}
Let $(M,g,f)$ be a complete connected noncompact non-Ricci-flat steady
gradient Ricci soliton of dimension $n\ge3$. The following conditions are equivalent:
\begin{enumerate}
\item[\textup{(i)}] $(M,g,f)$ is asymptotically cylindrical in the
sense of \cref{defn:asym_round}.
\item[\textup{(ii)}] There is $p_0\in M$ such that, as $t\to-\infty$,
\[
 (M,(-t)^{-1}g(t),p_0)\longrightarrow
 (\mathbb R\times\mathbb S^{n-1},\dd z^2+2(n-2)g_0,o)
\]
in the pointed smooth Cheeger--Gromov topology.
\item[\textup{(iii)}] The canonical Ricci flow arises as a noncollapsed
\(\mathbb F\)-limit of \(n\)-dimensional Ricci flows with complete time
slices and bounded curvature on compact time intervals. It admits a
tangent flow at infinity equal to the shrinking cylinder
\(\mathbb R\times\mathbb S^{n-1}\) in Bamler's \(\mathbb F\)-sense
\cite[Definition~6.55]{BamlerCompactness}, based at the conjugate heat flow
from a fixed spacetime point.
\end{enumerate}
Moreover, under these conditions, $g$ has bounded curvature, $f$ is proper,
and every sufficiently high level of $f$ is a spherical level set contained
in a cylindrical chart about the corresponding point of $\varphi_t(p_0)$.
\end{lemma}

\begin{proof}
\textup{(i)} $\Rightarrow$ \textup{(ii)}. Set $r=(n-1)/(2R_g)$.
The pointed smooth convergence (c.f. \cite{Brendle2013, Brendle2014, DengZhu2020a}) in \textup{(i)} gives, as $p\to\infty$,
\[
 \Delta R_g=o(R_g^2),\qquad \abs{\Ric_g}^2=\Bigl(\frac1{n-1}+o(1)\Bigr)R_g^2.
\]
By \eqref{eq:driftR},
\[
 \inner{\nabla f}{\nabla r}=\frac{n-1}{2R_g^2}\bigl(\Delta R_g+2\abs{\Ric_g}^2\bigr)\longrightarrow1.
\]
Let $K$ be a compact set outside which $\inner{\nabla f}{\nabla r}\ge1/2$.
Since $\abs{\nabla f}\le1$, the trajectory $\varphi_t(p)$ stays outside $K$
for $0\le t<d_g(p,K)$. Along it, $r$ decreases at rate at least $1/2$, so
$r(p)\ge d_g(p,K)/2$ by positivity of $r$. Thus $r(p)\to\infty$ as $p\to\infty$.
Choose $p_0$ with $r(p_0)$ large enough that
$\inner{\nabla f}{\nabla r}\ge1/2$ on $\{r\ge r(p_0)\}$.
Set $q(s):=\varphi_{-s}(p_0)$. For $r(p_0)$ sufficiently large,
$r(q(s))$ grows at least at rate $1/2$, so $q(s)\to\infty$.
Thus $(r\circ q)'(s)\to1$, and integration gives $r(q(s))/s\to1$.
Applying \textup{(i)} at $q(s)$ and using $g(-s)=\varphi_{-s}^*g$ implies \textup{(ii)}.

\smallskip
\textup{(ii)} $\Rightarrow$ \textup{(iii)}. Write $p_0(s)=\varphi_{-s}(p_0)$. Since
\[
 sR_g(p_0(s))\to\frac{n-1}{2},\qquad \frac{\dd}{\dd s}f(p_0(s))=1-R_g(p_0(s))\to1,
\]
$f(p_0(s))$ is strictly increasing for large $s$ and $f(p_0(s))/s\to1$.
Let $\Sigma_{f(p_0(s))}$ be the connected component of $\{f=f(p_0(s))\}$
containing $p_0(s)$, and set $u_s=(f-f(p_0(s)))/\sqrt{s}$. Then
\[
 \abs{\nabla_{s^{-1}g}u_s}_{s^{-1}g}^2=1-R_g\to1,\qquad \nabla_{s^{-1}g}^2u_s=s^{-1/2}\Ric_g\to0.
\]
Since $u_s(p_0(s))=0$ and $(M,s^{-1}g,p_0(s))$ converges to the cylinder
$(\mathbb R\times\mathbb S^{n-1},\dd z^2+2(n-2)g_0,o)$, with $z(o)=0$,
it follows that $u_s\to z$ locally in $C^1$.
In particular, $\Sigma_{f(p_0(s))}$ converges in $C^1$ to
$\{0\}\times\mathbb S^{n-1}$ in the pointed convergence charts.

For small $a$, the flow of $\nabla_g f/\abs{\nabla_g f}_g^2$ carries the
marked point along the curve $s\mapsto p_0(s)$ and maps $\Sigma_T$ diffeomorphically onto
$\Sigma_{T+a}$. In particular, the union of these spheres over any compact
interval of large levels is compact. Thus, for large $T_0$,
$U:=\bigcup_{T>T_0}\Sigma_T$ is open, connected and unbounded, with
$\partial U=\Sigma_{T_0}$. Since $\Sigma_{T_0}$ is a connected, two-sided
sphere and $M$ has one end by Munteanu--Wang~\cite[Corollary~1.1]{MunteanuWang2011},
$M\setminus U$ is compact. Hence $f$ is proper and $\Sigma_T=\{f=T\}$
for all sufficiently large $T$. Cylindrical convergence along $p_0(s)$ gives
$\sup_{\Sigma_{f(p_0(s))}}\abs{\operatorname{Rm}_g}_g=O(s^{-1})$, so $g$ has
bounded curvature.

Thus the canonical flow is a noncollapsed \(\mathbb F\)-limit, with itself
as the constant approximating sequence; see \cite[\S~2.1]{Bamlerstructure} and
\cite[Appendix~A]{BamlerFundamentalGroup}.
The distance condition following \cite[Definition~2.7]{BamlerLai} is automatic.
Then \cite[Lemma~6.38]{BamlerLai} gives \textup{(iii)}.

\smallskip
\textup{(iii)} $\Rightarrow$ \textup{(ii)}. This follows from \cite[Corollary~6.37]{BamlerLai}.

\smallskip
\textup{(ii)} $\Rightarrow$ \textup{(i)}. Given $p_i\to\infty$, choose $s_i\to\infty$ with $f(p_0(s_i))=f(p_i)$.
As shown in the proof of \textup{(ii)} $\Rightarrow$ \textup{(iii)}, the
level set $\Sigma_{f(p_0(s_i))}$ containing $p_i$ converges to
$\{0\}\times\mathbb S^{n-1}$ in the pointed convergence to the cylinder. Thus
$(M,s_i^{-1}g,p_i)$ converges locally smoothly to the cylinder.
In particular, since $s_iR_g(p_i)\to\frac{n-1}{2}$, it follows for $r_i=(n-1)/(2R_g(p_i))$ that $\frac{r_i}{s_i}\to1$.
This proves \textup{(i)}.
\end{proof}

By \cref{lem:asym-round-equivalence}, $f$ attains its minimum. We fix its
additive constant by
\begin{equation}\label{eq:f-normalization}
 \min_M f=0.
\end{equation}

\medskip
Finally, we recall some asymptotic properties of the Bryant soliton.
Let $(M_B,g_B,f_B)$ denote the complete normalized $n$-dimensional Bryant
steady gradient Ricci soliton, and let $o_B$ be its rotational center (also the unique critical point of $f_B$).  Thus
\[
 R_{g_B}+\abs{\nabla_{g_B}f_B}_{g_B}^2=1,
 \qquad f_B(o_B)=\min_{M_B}f_B=0.
\]
Let $r_B:[0,\infty)\longrightarrow(0,1]$ and $B:(0,1]\to(0,\infty)$ be functions defined as below,

so that for $T=f_B>0$,
\begin{equation}\label{eq:Bryant-profile-defs}
 R_{g_B}=r_B\circ f_B,
 \quad
 B(r_B(T))=\Ric_{g_B}(\nu_B,\nu_B)=-\frac12r_B'(T),
 \quad
 \nu_B=\frac{\nabla_{g_B}f_B}{\abs{\nabla_{g_B}f_B}_{g_B}}.
\end{equation}
Let $\rho_B=d_{g_B}(o_B,\cdot)$ and let
$\theta:M_B\setminus\{o_B\}\to\Sn^{n-1}$ be the angular projection.
In the product coordinates $(T,\theta)$, we have
\begin{equation}\label{eq:Bryant-metric}
 g_B=\frac{\dd T^2}{1-r_B(T)}+a_B^2(T)g_0,
\end{equation}
where $a_B:(0,\infty)\to(0,\infty)$ is the Bryant warping function.
Since $\dd\rho_B/\dd T=(1-r_B)^{-1/2}\to1$,  the asymptotic cylindricality implies as $T\to\infty$,
\begin{equation}\label{eq:Bryant-asymptotics}
 T\sim\rho_B,\qquad r_B(T)\sim\frac{n-1}{2T},\qquad a_B^2(T)\sim2(n-2)T.
\end{equation}
Here and below, $\sim$ means that the ratio tends to $1$ at the indicated limit.

By rotational symmetry, the restrictions of $f_B$ and $R_{g_B}$ to a
geodesic through $o_B$ are smooth and even in signed distance, and hence
are smooth functions of $\rho_B^2$ near the tip.
At $o_B$, $\Ric_{g_B}=g_B/n$ by rotational symmetry and $R_{g_B}(o_B)=1$,
so $f_B=\rho_B^2/(2n)+\Oterm(\rho_B^4)$. Thus $T=f_B$ has derivative
$1/(2n)$ with respect to $\rho_B^2$ at the tip. By the inverse function
theorem, $\rho_B^2$ is a smooth function of $T$ near zero. Composing with
the scalar curvature as a smooth function of $\rho_B^2$ shows that $r_B$ is smooth at $T=0$, with
\begin{equation}\label{eq:Bryant-tip}
 r_B(0)=1,\qquad r_B'(0)=-\frac{2}{n},
\end{equation}
where the derivative follows by letting $T\to0$ in the second equation in \eqref{eq:Bryant-profile-defs}.
Since $r_B'(0)\ne0$, the inverse function theorem implies that
$r_B^{-1}$ is smooth at $s=1$. Therefore
$B=-\frac12(r_B'\circ r_B^{-1})$ is smooth there.

The Bryant asymptotics~\cite[Appendix~A]{BDNS} give
$\Ric_{g_B}(\nu_B,\nu_B)\sim(n-1)/(4\rho_B^2)$ as $\rho_B\to\infty$.
Together with \eqref{eq:Bryant-profile-defs} and $T\sim\rho_B$,
this gives
\begin{equation}\label{eq:Bryant-rB-derivative-asymptotic}
 -r_B'(T)\sim\frac{n-1}{2T^2}\qquad\text{as }T\to\infty,
\end{equation}
and, as $s\downarrow0$,
\begin{equation}\label{eq:B-endpoints}
 B(s)\sim\frac{s^2}{n-1}.
\end{equation}

\smallskip
The following lemma records the ODE satisfied by $r_B$, which follows from \eqref{eq:driftR}.

\begin{lemma}\label{lem:Bryant-model}
Let $r_B$, $B$, and $a_B$ be functions defined above. Then
\begin{enumerate}
\item[\textup{(i)}] For $T>0$ and $s=r_B(T)$,
\begin{equation}\label{eq:Bryant-rB}
 (1-s)r_B''(T)+r_B'(T)=-2\left(B^2(s)+\frac{(s-B(s))^2}{n-1}\right).
\end{equation}

\item[\textup{(ii)}] For every $j\ge0$, $r_B^{(j)}(T)$ and $a_B^{(j)}(T)$
are bounded by polynomial functions in $T$ for sufficiently large $T$.
\end{enumerate}
\end{lemma}

\begin{proof}
For \textup{(i)}, the tangential Ricci eigenvalue is $(s-B(s))/(n-1)$ by the definition of $B$. Applying
\eqref{eq:driftR} to $R_{g_B}=r_B\circ f_B$ with
$\Delta_{f_B}f_B=1$ and $\abs{\nabla_{g_B}f_B}_{g_B}^2=1-s$ gives
\eqref{eq:Bryant-rB}. Moreover, substituting $r_B'(T)=-2B(s)$ and
$r_B''(T)=4B(s)B'(s)$ in \eqref{eq:Bryant-rB} gives
\begin{equation}\label{eq:Bryant-B}
 2(1-s)B(s)B'(s)=B(s)-B^2(s)-\frac{(s-B(s))^2}{n-1}.
\end{equation}

For \textup{(ii)}, differentiating \eqref{eq:Bryant-B} $j$ times gives
$B^{(j+1)}$ with coefficient $2(1-s)B(s)$, while all other terms involve
derivatives of order at most $j$; since $B(s)\sim s^2/(n-1)$, division by
this coefficient introduces an $O(s^{-2})$ factor, so by induction
each derivative of $B$ is bounded by a power of $s^{-1}$ near $s=0$.
Since $r_B(T)\sim(n-1)/(2T)$, the derivative bounds for $B$ give polynomial
bounds for all derivatives of $r_B$ by differentiating $r_B'=-2B(r_B)$
from \eqref{eq:Bryant-profile-defs} by induction.
In the coordinates of \eqref{eq:Bryant-metric}, $\nabla_{g_B}f_B=(1-r_B)\partial_T$,
so for a $g_B$-unit vector $X$ tangent to $\{f_B=T\}$,
\[
 \nabla_{g_B}^2f_B(X,X)=\tfrac12(\mathcal L_{\nabla_{g_B}f_B}g_B)(X,X)=\frac{1-r_B}{2}\partial_T(a_B^2g_0)(X,X)=(1-r_B)\frac{a_B'}{a_B}.
\]
Equating this to $\Ric_{g_B}(X,X)=(r_B-B(r_B))/(n-1)$ yields polynomial
bounds for all derivatives of $a_B$ by induction,
using $1-r_B\to1$, $a_B=O(T^{1/2})$, and the polynomial bounds for the derivatives of $r_B,B$.
\end{proof}

\section{Exponential convergence estimate}\label{s:exp}
In this section, we apply Bamler--Lai's theorem on exponential Bryant
asymptotics~\cite[Theorems~1.3 and~7.1]{BamlerLai} to compare the metrics
on the ends of \(M\) and \(M_B\). We then use the soliton identities to
compare their potentials.

Fix \(p=p_0\) as in \cref{lem:asym-round-equivalence}. The pointed metrics
\((M,(-t)^{-1}g(t),p)\) converge to the cylinder as \(t\to-\infty\).

By \cref{lem:asym-round-equivalence} the canonical flow has complete time
slices and uniformly bounded curvature. By \cite[Proposition~2.1]{Lu2010},
the canonical flow satisfies the pseudolocality property of
\cite[Definition~2.4]{BamlerLai}, with a smaller dimensional constant if necessary.

Let \(\mu_1(t)\le\cdots\le\mu_n(t)\) be the Ricci eigenvalues of \(g\)
at \(\varphi_t(p)\), equivalently those of \(g(t)\) at \(p\).
Smooth cylindrical convergence gives \((-t)\mu_1(t)\to0\) and
\((-t)\mu_j(t)\to\tfrac12\) for \(j\ge2\), so \(\abs{\mu_1(t)}=\min_i\abs{\mu_i(t)}\)
for sufficiently negative \(t\). By \eqref{eq:gradR}, \eqref{eq:normalization},
and the constancy of scalar curvature on the limiting cylinder,
\[
 \abs{\mu_1(t)}
 \le \frac{\abs{\nabla R_g}_g(\varphi_t(p))}{2\abs{\nabla f}_g(\varphi_t(p))}
 =o((-t)^{-3/2}).
\]
Thus, along the cylindrical worldline $t\mapsto(p,t)$,
\[
 (-t)\log(-t)\abs{\mu_1(t)}\longrightarrow0\qquad(t\to-\infty),
\]
which is the $k=1$ condition in Bamler--Lai's Ricci-eigenvalue
criterion~\cite[Lemma~6.65]{BamlerLai}. Since $\sup_M R_{g(t)}$ is independent
of $t$, the flow is not homothetic to the shrinking cylinder. Hence
it has dominant linear mode in the sense of \cite[Definition~6.64]{BamlerLai}.

\begin{lemma}
\label{lem:potential-alignment}
Under the normalization \eqref{eq:normalization}, there
exist compact domains \(K\subset M\) and \(K_B\subset M_B\), a
diffeomorphism $$\phi:M\setminus K\to M_B\setminus K_B,$$ and constants
\(c_\infty\in\mathbb R\) and \(\beta_1>0\) such that, identifying these ends
via \(\phi\), we have that
\begin{align}
 \abs{\nabla_{g_B}^j(g-g_B)}_{g_B}
 &\le C_j\e^{-\beta_1\rho_B}\qquad\text{for every }j\ge0,\notag\\
 \abs{\nabla_{g_B}^j(f_B-f-c_\infty)}_{g_B}
 &\le C_j\e^{-\beta_1\rho_B}\qquad\text{for every }j\ge0,
 \label{eq:potential-alignment}\\
 \big|R_g-r_B(f+c_\infty)\big|&\le C_1\e^{-\beta_1\rho_B},
 \label{eq:R-alignment}
\end{align}
holds everywhere on $M_B\setminus K_B$.
\end{lemma}

\begin{proof}
By the exponential estimate theorem in Bamler-Lai \cite[Theorem~7.1]{BamlerLai}, with
\(k=1\), there exist \(\Lambda>0\)
and a time-zero diffeomorphism onto its image \(\psi:\{\rho_B>Y_0\}\to M\)
such that \(h=\psi^*(\Lambda g)\) converges exponentially to \(g_B\) in every
derivative order, with rates that may depend on the order.
Interpolation on unit \(g_B\)-balls between the \(C^0\) exponential
bound with rate \(c_0\) and the uniform covariant \(C^{2j}\) bounds gives
\begin{equation}\label{eq:metric-exponential-input}
 \abs{\nabla_{g_B}^j(h-g_B)}_{g_B}\le C_j\e^{-\beta\rho_B}
\end{equation}
for every \(j\ge0\) on the far end, with \(\beta=c_0/2\) independent of \(j\).
By \eqref{eq:metric-exponential-input}, $R_g\circ\psi=\Lambda R_h\to0$
as $\rho_B\to\infty$. Since $R_g>0$, the image escapes compact sets as
$\rho_B\to\infty$. Thus, for large $Y_1>Y_0$, $\psi(\{\rho_B>Y_1\})$ is an
unbounded component of the complement of $\psi(\{\rho_B=Y_1\})$. By the
one-endedness used in \cref{lem:asym-round-equivalence}, the domains
\(K_B:=\{\rho_B\le Y_1\}\) and \(K:=M\setminus\psi(\{\rho_B>Y_1\})\) are compact.
We may assume $\psi=\textnormal{id}$ in the following.
Then \(h=\Lambda g\) on $M_B\setminus K_B$.

By \eqref{eq:B-endpoints}, \eqref{eq:Bryant-asymptotics}, and
\eqref{eq:metric-exponential-input}, both \(\Ric_{g_B}\) and \(\Ric_h\)
are bounded below by \(c(1+\rho_B)^{-2}g_B\) for large \(\rho_B\).
Constant rescaling preserves the soliton equation, so applying and differentiating
\eqref{eq:gradR} to \((h,f)\) and \((g_B,f_B)\), together with the
exponential curvature estimates from \eqref{eq:metric-exponential-input}
and the polynomial curvature derivative bounds, gives
\[
 \abs{\nabla_{g_B}^{j+1}(f-f_B)}_{g_B}\le C_j\e^{-\beta_1\rho_B}\qquad(j\ge0),
\]
with \(\beta_1=\beta/2\).
In particular, by \eqref{eq:normalization},
\[
 \Lambda^{-1}=R_h+\abs{\nabla_h f}_h^2
 \longrightarrow 1,
\]
where the limit follows by comparison with
\(R_{g_B}+\abs{\nabla_{g_B}f_B}_{g_B}^2=1\). Hence \(\Lambda=1\), so
\eqref{eq:metric-exponential-input} holds with \(h=g\).

Integrating the estimate for \(\nabla_{g_B}(f_B-f)\) along Bryant radial rays gives
a limit on each ray, with a uniform error \(\Oterm(\e^{-\beta_1\rho_B})\).
The oscillation on \(\{f_B=T\}\) is at most
\(C a_B(T)\e^{-\beta_1\rho_B(T)}\to0\), since the level has intrinsic diameter
\(\pi a_B(T)\). Thus the limits agree and define \(c_\infty\).
Together with the estimates for \(\nabla_{g_B}(f_B-f)\), this proves
\eqref{eq:potential-alignment}. Finally, the metric estimate gives
\(R_g=r_B(f_B)+\Oterm(\e^{-\beta\rho_B})\).
Since \(r_B'\) is bounded, \eqref{eq:R-alignment} follows from
\eqref{eq:potential-alignment} and the mean-value theorem.
Enlarge \(K_B\) and \(K\) correspondingly so that \(f+c_\infty>0\)
on the end. Both estimates then hold throughout \(M_B\setminus K_B\)
after increasing the constants on any remaining compact annuli.
\end{proof}

\section{Scalar curvature comparison}\label{s:comparison}
Using a maximum principle argument, we compare $R_g$ with $r_B(f+c)$ by decreasing $c$ from a sufficiently large value to $c_\infty$, determined by the asymptotics established in the previous section. The shifting of $c$ is partly motivated by the extrinsic sliding argument in the mean curvature flow (see for instances \cite{BamlerLaiMCF1}). The resulting comparison will be used in the next section to show that, after a suitable conformal change of $g$, the metric still has nonnegative scalar curvature, which is crucial for applying the positive mass theorem.

\begin{lemma}\label{lem:sliding-identity}
For $c\ge0$, put $Z_c=R_g-r_B(f+c)$ and $s=r_B(f+c)$. Set
\[
 \lambda=\frac{s-B(s)}{n-1},\qquad k=\frac{B(s)-\lambda}{1-s},
 \qquad S_c=\lambda g+k\,\dd f\otimes\dd f.
\]
With $k$ extended smoothly at $s=1$, we have on $M$
\begin{equation}\label{eq:sliding-identity}
\begin{aligned}
 \Delta_f Z_c-2k\inner{\nabla f}{\nabla Z_c}
 &=\bigl(r_B''(f+c)-4\lambda+2k^2Z_c\bigr)Z_c
   -2\abs{\Ric_g-S_c}^2.
\end{aligned}
\end{equation}
\end{lemma}
\begin{Remark} 
The tensor \(S_c\) is modeled on the Bryant Ricci tensor. In the Bryant formula \(\Ric_{g_B}=\lambda g_B+k\,df_B\otimes df_B\), we replace \(g_B\) and \(df_B\) by \(g\) and \(df\), and evaluate \(\lambda\) and \(k\) at the matched level \(f_B=f+c\).
\end{Remark}
\begin{proof}
By rotational symmetry at $o_B$, $B(1)=\lambda(1)$, so 
$B(s)-\lambda$ vanishes at $s=1$. Since $B(s)-\lambda$ is also smooth at $s=1$ as shown after \eqref{eq:Bryant-tip}, $k$ extends smoothly there.

By \eqref{eq:normalization} and \eqref{eq:trace-soliton}, $\Delta_f f=1$.
Thus \eqref{eq:driftR} and the scalar curvature equation in
\cref{lem:Bryant-model}, \eqref{eq:Bryant-rB}, give
\begin{align}
 \Delta_f Z_c
 &=-2\abs{\Ric_g}^2-(1-R_g)r_B''(f+c)-r_B'(f+c)\notag\\
 &=r_B''(f+c)Z_c-2\bigl(\abs{\Ric_g}^2-B^2(s)-(n-1)\lambda^2\bigr).
 \label{eq:comparison-drift-Z}
\end{align}
By \eqref{eq:gradR} and $r_B'(f+c)=-2B(s)$,
\begin{equation}\label{eq:comparison-radial-ricci}
 \Ric_g(\nabla f,\nabla f)
 =B(s)(1-R_g)-\tfrac12\inner{\nabla f}{\nabla Z_c}.
\end{equation}

We first subtract $\lambda g$ from the Ricci tensor and make the analogous
subtraction on Bryant, where only the radial eigenvalue $B(s)-\lambda$
remains. Using $\operatorname{tr}_g\Ric_g=R_g$ and $s=B(s)+(n-1)\lambda$, we have
\begin{align*}
 \abs{\Ric_g}^2-\abs{\Ric_g-\lambda g}^2
 &=2\lambda R_g-n\lambda^2,\\
 \bigl[B^2(s)+(n-1)\lambda^2\bigr]-(B(s)-\lambda)^2
 &=2\lambda s-n\lambda^2.
\end{align*}
Subtracting the second line from the first and using $R_g-s=Z_c$ gives
\begin{equation}\label{eq:comparison-tangential-norm}
 \abs{\Ric_g}^2-B^2(s)-(n-1)\lambda^2
 =\abs{\Ric_g-\lambda g}^2-(B(s)-\lambda)^2+2\lambda Z_c.
\end{equation}

We next subtract the radial term $k\,\dd f\otimes\dd f$. Expanding the square gives
\[
 \abs{\Ric_g-S_c}^2
 =\abs{\Ric_g-\lambda g}^2-2k(\Ric_g-\lambda g)(\nabla f,\nabla f)+k^2\abs{\nabla f}^4.
\]
By \eqref{eq:comparison-radial-ricci} and $1-R_g=\abs{\nabla f}^2$, the cross term is
\[
 -2k(\Ric_g-\lambda g)(\nabla f,\nabla f)
 =-2k(B(s)-\lambda)\abs{\nabla f}^2+k\inner{\nabla f}{\nabla Z_c}.
\]
Completing the square then gives
\begin{align}
 \abs{\Ric_g-S_c}^2
 &=\abs{\Ric_g-\lambda g}^2+k\inner{\nabla f}{\nabla Z_c}
   +\bigl(B(s)-\lambda-k\abs{\nabla f}^2\bigr)^2-(B(s)-\lambda)^2\notag\\
 &=\abs{\Ric_g-\lambda g}^2-(B(s)-\lambda)^2
   +k\inner{\nabla f}{\nabla Z_c}+k^2 Z_c^2,\label{eq:comparison-radial-norm}
\end{align}
where the last line uses $B(s)-\lambda-k\abs{\nabla f}^2
=k(1-s)-k(1-R_g)=kZ_c$. Combining the two norm comparisons \eqref{eq:comparison-tangential-norm}
and \eqref{eq:comparison-radial-norm} in \eqref{eq:comparison-drift-Z} gives
\eqref{eq:sliding-identity}.
\end{proof}

\begin{proposition}\label{prop:scalar-comparison}
We have $c_\infty\ge0$ and
\begin{equation}\label{eq:R-comparison}
 R_g\ge r_B(f+c_\infty)\qquad\text{on }M.
\end{equation}
\end{proposition}

\begin{proof}
For $c\ge0$, use the notation of \cref{lem:sliding-identity}.
Suppose $Z_c\ge0$ on $M$. Near any zero of $Z_c$, by \eqref{eq:sliding-identity} there exists $C\ge0$
such that
$\Delta_f Z_c-2k\inner{\nabla f}{\nabla Z_c}-CZ_c\le0$. Since the zeroth-order
coefficient is $-C\le0$, the strong minimum principle~\cite[Theorem~3.5]{GilbargTrudinger}
gives $Z_c\equiv0$ in a connected neighborhood of each zero. Hence the zero
set is open. Since $M$ is connected,
either $Z_c>0$ everywhere or $Z_c\equiv0$.

By \eqref{eq:R-alignment} and $f\sim\rho_B$,
$R_g-r_B(f+c_\infty)=o(f^{-2})$. The mean-value theorem and
\eqref{eq:Bryant-rB-derivative-asymptotic} then give, for every fixed $c\ge0$,
\begin{equation}\label{eq:sliding-end-sign_0}
 f^2Z_c\longrightarrow\frac{n-1}{2}(c-c_\infty)\qquad\text{as }f\to\infty.
\end{equation}

Fix $c_0>\max\{0,c_\infty\}$ and choose a compact set $K$ such that
$Z_{c_0}>0$ on $M\setminus K$ by \eqref{eq:sliding-end-sign_0}.
Since $r_B$ is decreasing, $Z_c\ge Z_{c_0}>0$ on $M\setminus K$ for every
$c\ge c_0$. The set of $c\in[c_0,\infty)$ for which $Z_c\ge0$ on $M$ is
closed by continuity and nonempty since $\min_K R_g>0$ and
$r_B(f+c)\le r_B(c)\to0$ as $c\to\infty$.
For each such $c$, the strong minimum principle gives $Z_c>0$ on $M$.
Its positive minimum on $K$ persists under small changes in $c$, so the
set is also relatively open. By connectedness of $[c_0,\infty)$,
$Z_{c_0}>0$ on $M$.

Letting $c_0\downarrow\max\{0,c_\infty\}$ gives
$Z_{\max\{0,c_\infty\}}\ge0$. If $c_\infty<0$, then $Z_0$ vanishes at a
minimum of $f$ by \eqref{eq:normalization} and \eqref{eq:Bryant-tip}, but
is positive near infinity by \eqref{eq:sliding-end-sign_0}. This contradicts
the strong minimum principle. Hence $c_\infty\ge0$ and $Z_{c_\infty}\ge0$.
\end{proof}

\section{Conformal deformation and the positive mass theorem}\label{s:conformal}
In this section, we complete the proof of Theorem \ref{thm:main}. Using the diffeomorphism constructed in Lemma \ref{lem:potential-alignment}, we find conformal factors that transform $g_B$ into the Euclidean metric $g_{\mathrm{Euc}}$ and $g$ into an asymptotically flat metric $\widetilde g$ of nonnegative scalar curvature, respectively. The exponential Bryant asymptotic of $g$ will imply that $\widetilde g$ converge to $g_{\mathrm{Euc}}$ sufficiently fast near infinity, thereby proving the vanishing of the ADM mass of $\widetilde g$. Our main theorem will follow from the rigidity of the positive mass theorem and the classification of locally conformally flat steady Ricci soliton. 

We first construct the conformal factor that maps the Bryant metric to Euclidean.

\begin{lemma}\label{lem:Bryant-conformal-factor}
There is a smooth positive function $u_B:[0,\infty)\to(0,\infty)$ with
$u_B(0)=1$ such that
\begin{equation}\label{eq:uB-def}
 \varrho(T):=a_B(T)u_B^{\frac{2}{n-2}}(T)
\end{equation}
is the Euclidean radius for the conformally changed metric. In the
coordinates of \eqref{eq:Bryant-metric},
\begin{equation}\label{eq:Bryant-flat}
 (u_B^{\frac{4}{n-2}}\circ f_B)g_B
 =\dd(\varrho\circ f_B)^2+(\varrho\circ f_B)^2g_0,
\end{equation}
which is isometric to the Euclidean metric $g_{\mathrm{Euc}}$ and
\begin{equation}\label{eq:rho-def}
 \frac{\dd}{\dd T}\log\varrho(T)=\frac{1}{a_B(T)\sqrt{1-r_B(T)}} \qquad (T>0).
\end{equation}
\end{lemma}

\begin{proof}
Matching the spherical coefficients in \eqref{eq:Bryant-flat} gives
\eqref{eq:uB-def}. Matching the radial coefficients with $\varrho'>0$ gives
$\varrho'=u_B^{\frac{2}{n-2}}/\sqrt{1-r_B}$, so
\begin{equation}\label{eq:lnub-def}
 \frac{\dd}{\dd T}\log u_B
 =\frac{n-2}{2}\left(\frac{1}{a_B\sqrt{1-r_B}}-\frac{a_B'}{a_B}\right).
\end{equation}
Near the tip, $a_B/\rho_B$ and $T/\rho_B^2$ are smooth functions of $\rho_B^2$
with values $1$ and $1/(2n)$ at the tip. Since $T$ is a smooth coordinate
in $\rho_B^2$ there, $a_B(T)/\sqrt{2nT}$ is smooth at $0$ with value $1$.
The same holds for $n(1-r_B(T))/(2T)$ by \eqref{eq:Bryant-tip}.
It follows that $(a_B\sqrt{1-r_B})^{-1}-1/(2T)$ is smooth at $0$, as is
$a_B'/a_B-1/(2T)=(\log(a_B/\sqrt{2nT}))'$. Thus the right-hand side of \eqref{eq:lnub-def} is
smooth at $0$, and integration with $u_B(0)=1$ gives a smooth positive
$u_B$ on $[0,\infty)$.

With $\varrho$ given by \eqref{eq:uB-def}, the equation \eqref{eq:lnub-def} for $u_B$ gives
\eqref{eq:rho-def}, and the coefficients agree in \eqref{eq:Bryant-flat}.
By \eqref{eq:rho-def} and \eqref{eq:Bryant-asymptotics}, $\varrho$ increases
from $0$ to $\infty$. Since $\varrho/\rho_B$ is positive and smooth in
$\rho_B^2$ at the tip, the Euclidean coordinates extend across $o_B$.
\end{proof}

For either pair $(h,F)=(g_B,f_B)$ or $(h,F)=(g,f+c_\infty)$, put $u=u_B(F)$.
By \eqref{eq:trace-soliton} and \eqref{eq:normalization},
$\Delta_h u=(1-R_h)u_B''(F)+R_hu_B'(F)$. The conformal scalar-curvature
formula therefore gives, with $c_n=\frac{4(n-1)}{n-2}$,
\begin{equation}\label{eq:conformal-soliton}
 u^{\frac{n+2}{n-2}}R_{u^{\frac{4}{n-2}}h}
 =-c_n(1-R_h)u_B''(F)+R_h\bigl[u_B(F)-c_nu_B'(F)\bigr].
\end{equation}
For $(h,F)=(g_B,f_B)$, the left-hand side vanishes by \eqref{eq:Bryant-flat}.
Using $R_{g_B}=r_B\circ f_B$ gives
\begin{equation}\label{eq:uB-ode}
 -c_n(1-r_B)u_B''-c_nr_Bu_B'+r_Bu_B=0.
\end{equation}
Evaluating at $T=0$ and using \eqref{eq:Bryant-tip} and $u_B(0)=1$ gives $u_B'(0)=1/c_n$.

The following lemma establishes the positivity needed in the scalar-curvature
factorization below.

\begin{lemma}\label{lem:Qpositive}
We have $u_B\ge1$, $u_B''>0$, and
\begin{equation}\label{eq:Qn}
 Q_n:=c_nu_B''-c_nu_B'+u_B
 =\frac{c_nu_B''}{r_B}>0.
\end{equation}
\end{lemma}

\begin{proof}
Rewrite \eqref{eq:uB-ode} as
\[
 (1-r_B)u_B''=r_B\left(\frac{u_B}{c_n}-u_B'\right).
\]
Differentiating the display at $T=0$ and using \eqref{eq:Bryant-tip}
and $u_B(0)/c_n=u_B'(0)=1/c_n$ gives
\[
 \left(\frac{u_B}{c_n}-u_B'\right)'(0)=\frac{2}{n}u_B''(0)=\frac{1}{c_n^2}-u_B''(0).
\]
Solving gives
\[
 \left(\frac{u_B}{c_n}-u_B'\right)'(0)=\frac{2}{(n+2)c_n^2}>0.
\]
Thus $u_B''(0)>0$ and $u_B'<u_B/c_n$ for small $T>0$.
At a first point of equality with $T>0$, \eqref{eq:uB-ode} gives $u_B''=0$, so
$(u_B/c_n-u_B')'=u_B/c_n^2>0$. This contradicts the derivative being
nonpositive at the first zero. Hence $u_B'<u_B/c_n$ for all $T>0$, so
$u_B''>0$ there by the same equation. Since $u_B'(0)>0$, we have $u_B\ge u_B(0)=1$.
Finally, \eqref{eq:uB-ode} gives $r_B Q_n=c_nu_B''$, proving \eqref{eq:Qn}.
\end{proof}

On $(M,g,f)$, set
\begin{equation}\label{eq:conformal-unknown}
 u=u_B(f+c_\infty),
 \qquad
 \widetilde g=u^{\frac{4}{n-2}}g.
\end{equation}
Subtracting \eqref{eq:uB-ode} at $T=f+c_\infty$ from
\eqref{eq:conformal-soliton} for $(h,F)=(g,f+c_\infty)$ gives
\begin{equation}\label{eq:conformal-factorization}
 u^{\frac{n+2}{n-2}}R_{\widetilde g}
 =\bigl[R_g-r_B(f+c_\infty)\bigr]Q_n(f+c_\infty).
\end{equation}
By \eqref{eq:R-comparison} and \cref{lem:Qpositive},
\begin{equation}\label{eq:Rtilde-positive}
 R_{\widetilde g}\ge0 \qquad\text{at every point of }M.
\end{equation}

The following lemma constructs an asymptotically flat end chart in which
$\widetilde g$ converges to $g_{\mathrm{Euc}}$ to every algebraic order.

\begin{lemma}\label{lem:AF}
Use the end identification of \cref{lem:potential-alignment} and the Bryant
coordinates \((T,\theta)\), with \(T=f_B\). For sufficiently large \(T_0\),
\[
 x:\{T>T_0\}\longrightarrow
 \mathbb R^n\setminus\overline{B_{\varrho(T_0)}(0)},
 \qquad x(T,\theta)=\varrho(T)\theta
\]
is a diffeomorphism. Its domain has compact complement in \(M\), and
\begin{equation}\label{eq:AF-decay}
 \abs{\partial_x^j\bigl((x^{-1})^*\widetilde g-g_{\mathrm{Euc}}\bigr)}_{g_{\mathrm{Euc}}}
 \le C_{j,N}\abs{x}^{-N-j}
 =C_{j,N}\varrho^{-N-j}(T)
\end{equation}
for every fixed $j,N$.
Consequently,
\begin{equation}\label{eq:AF-consequences}
 R_{\widetilde g}\in L^1(M,\dd V_{\widetilde g}),
 \qquad
 m_{\mathrm{ADM}}(\widetilde g)=0.
\end{equation}
\end{lemma}

\begin{proof}
Set $\lambda_B(T):=u_B^{\frac{4}{n-2}}(T)$, and $\widetilde g_B:=\lambda_B(T)g_B=x^*g_{\mathrm{Euc}}$,
where the last equality follows from \eqref{eq:Bryant-flat}.
By \eqref{eq:rho-def} and \eqref{eq:uB-def},
\[
 (\log\lambda_B)'=\frac{2}{a_B\sqrt{1-r_B}}-2\frac{a_B'}{a_B}.
\]
Together with the derivative bounds in \cref{lem:Bryant-model}, this gives
polynomial bounds for all derivatives of \(\dd\log\lambda_B\).
Put $a:=\frac{\lambda_B(f+c_\infty)}{\lambda_B(T)}$.
Then
\[
 \log a=(f+c_\infty-T)\int_0^1
 (\log\lambda_B)'\bigl(T+s(f+c_\infty-T)\bigr)\,\dd s.
\]
So by the potential exponential estimate in \eqref{eq:potential-alignment}, and the polynomial derivative
bounds of $\log\lambda_B$, there is a fixed \(\gamma>0\) such that for all $j\ge 0$,
\begin{equation*} 
\abs{\nabla_{g_B}^j(a-1)}_{g_B}\le C_j\e^{-\gamma T}.
\end{equation*}
Define $E:=\lambda_B(T)^{-1}\bigl(\widetilde g-\widetilde g_B\bigr)
 =a(g-g_B)+(a-1)g_B$.
Then the estimate of $a-1$ and the exponential estimate of $g-g_B$ in \eqref{eq:metric-exponential-input} give
\begin{equation}\label{eq:conformal-end-error}
 \abs{\nabla_{g_B}^j E}_{g_B}\le C_j\e^{-\gamma T}\qquad(j\ge0).
\end{equation}

Recall that \(\widetilde g_B=\lambda_Bg_B=x^*g_{\mathrm{Euc}}\). To convert Bryant
derivatives to Euclidean derivatives, we use the conformal connection formula
for arbitrary smooth vector fields \(U,V\) on the Bryant end:
\[
 (\nabla_{\widetilde g_B}-\nabla_{g_B})_U V
 =\frac12\bigl(\dd\log\lambda_B (U)\, V+\dd\log\lambda_B (V)\, U-g_B(U,V)\nabla_{g_B}\log\lambda_B\bigr).
\]
Since \(\nabla_{g_B}g_B=0\), for each \(j\ge0\) there are constants \(C_j,d_j\) such that
\[
 \abs{\nabla_{g_B}^{j}(\nabla_{\widetilde g_B}-\nabla_{g_B})}_{g_B}
 \le C(n)\abs{\nabla_{g_B}^{j+1}\log\lambda_B}_{g_B}
 \le C_j(1+T)^{d_j}.
\]
Iterating the connection-difference formula
\cite[(2.48)]{ChowLuNi} and applying the Leibniz rule gives the following formula that relates Riemanninan connections $\nabla_{g_0}$ and $\nabla_{g_1}$ of any two metrics $g_0$ and $g_1$: for any smooth covariant tensor $\alpha$
and integer $k\geq 1$,
\[
\left(\nabla^k_{g_1}-\nabla_{g_0}^k\right)\alpha=\sum_{i=0}^{k-1}\sum_{l=1}^{k-i}\sum_{i+l+j_1+\cdots+j_l=k}\nabla_{g_1}^i \alpha\ast\nabla_{g_1}^{j_1}A\ast\cdots\ast\nabla_{g_1}^{j_l}A,
\]
where $A:=\nabla_{g_1}-\nabla_{g_0}$. Substitute $g_1=g_B$, $g_0=\widetilde g_B$, $\alpha=\lambda_B(T)E,$ and
set \(S_j:=\lambda_B^{-1}\nabla_{\widetilde g_B}^j(\lambda_B E)\).

We can express \(S_j\)
in terms of \(\nabla_{g_B}^\ell E\), \(\ell\le j\), with coefficients formed
from the connection difference, \(\dd\log\lambda_B\), and their
\(g_B\)-covariant derivatives. These coefficients are polynomially bounded
in \(T\) by the preceding estimates.
By \eqref{eq:conformal-end-error}, \(\abs{S_j}_{g_B}\) is bounded by a polynomial
in \(T\) times \(\e^{-\gamma T}\).
Since \(S_j\) has \(j+2\) covariant indices,
\[
 \abs{\nabla_{\widetilde g_B}^j(\lambda_B E)}_{\widetilde g_B}
 =\lambda_B^{-j/2}\abs{S_j}_{g_B}
 =\left(\frac{a_B}{\varrho}\right)^j\abs{S_j}_{g_B},
\]
where the last equality uses \(\lambda_B=(\varrho/a_B)^2\) by \eqref{eq:uB-def}.
The map \(x\) takes \(\nabla_{\widetilde g_B}\) to Cartesian Euclidean
differentiation, and \(a_B^j=\Oterm(T^{j/2})\).
Thus, for each \(j\), there are \(d_j\ge0\) and \(C_j\) such that
\[
 \abs{\partial_x^j\bigl((x^{-1})^*\widetilde g-g_{\mathrm{Euc}}\bigr)}_{g_{\mathrm{Euc}}}
 \le C_j\varrho^{-j}(T)(1+T)^{d_j}\e^{-\gamma T}.
\]
By \eqref{eq:rho-def} and \eqref{eq:Bryant-asymptotics},
\begin{equation}\label{eq:rho-asymptotic}
 \log\varrho(T)=\sqrt{\frac{2T}{n-2}}+o(\sqrt T).
\end{equation}
Hence
\begin{equation}\label{eq:exp-absorbs-polynomial}
 (1+T)^d\e^{-\gamma T}=\Oterm(\varrho^{-N}(T))
 \qquad\text{for every fixed }N,d\ge0.
\end{equation}
This proves \eqref{eq:AF-decay}. Taking $N>n-2$, the estimates with $j\le2$
imply that $R_{\widetilde g}\in L^1(M,\dd V_{\widetilde g})$, while the
estimate with $j=1$ makes the ADM flux integral~\cite{ADM59, Bartnik} tend to zero, i.e.
\[
\frac{1}{\rho}\int_{\mathbb S^{n-1}(\rho)}\left(\partial_{x^i}\left((x^{-1})^*\widetilde g_{ij}\right)-\partial_{x^j}\left((x^{-1})^*\widetilde g_{ii}\right)\right)x^j\, d\sigma \to 0  \text{  as  } \rho\to\infty,
\]
where $d\sigma$ is the surface measure on the Euclidean sphere $\mathbb S^{n-1}(\rho)$ of radius $\rho$.
This proves \eqref{eq:AF-consequences}.
\end{proof}

\begin{proof}[Proof of Theorem \ref{thm:main}:]
The metric $\widetilde g$ is complete because
\[
 \widetilde g=u^{\frac{4}{n-2}}g\ge g,
\]
where $u=u_B(f+c_\infty)\ge1$ and $g$ is complete.
By properness of $f$ and \cref{lem:AF}, $(M,\widetilde g)$ has a single
asymptotically flat end, with $R_{\widetilde g}\in L^1(M,\dd V_{\widetilde g})$
and $m_{\mathrm{ADM}}(\widetilde g)=0$. Since $R_{\widetilde g}\ge0$ by
\eqref{eq:Rtilde-positive}, the smooth case of the zero-mass rigidity
theorem~\cite[Theorem~1.1]{KhuriWangWang} implies that $(M,\widetilde g)$ is
isometric to Euclidean space.

Since $\widetilde g=u^{\frac{4}{n-2}}g$, the metric $g$ is locally conformally flat.
Since $(M,g)$ is complete, noncompact, and nonflat, the classification of
Cao--Chen~\cite[Theorem~1.2]{CaoChen} (see also Catino--Mantegazza~\cite[Proposition~2.3]{CatinoMantegazza}) implies that
$(M,g)$ is isometric to the Bryant soliton up to scaling.
This proves \cref{thm:main}.
\end{proof}

\end{document}